\documentclass{article}

\usepackage{amsfonts,amsmath,amssymb,amsthm}
\usepackage{natbib}
\usepackage[dvipsnames]{xcolor}
\usepackage{titletoc}
\usepackage{parskip}
\definecolor{color1}{HTML}{105e8a}
\definecolor{color3}{HTML}{b82a0c}
\definecolor{color5}{HTML}{80037e}
\usepackage{hyperref}
\hypersetup{
    colorlinks,
    linkcolor={color1},
    citecolor={color3},
    urlcolor={color5}
}
\usepackage{authblk}
\usepackage{microtype}
\usepackage{tikz}

\usepackage[noabbrev, capitalize]{cleveref}

\newtheoremstyle{thmstyle}
  {6pt} 
  {2pt} 
  {\itshape} 
  {} 
  {\bfseries} 
  {.} 
  {.5em} 
  {} 

\newtheoremstyle{defstyle}
  {6pt} 
  {2pt} 
  {} 
  {} 
  {\bfseries} 
  {.} 
  {.5em} 
  {} 

\theoremstyle{thmstyle}
\newtheorem{theorem}{Theorem}
\crefname{theorem}{Theorem}{Theorems}
\newtheorem{lemma}{Lemma}
\newtheorem{proposition}{Proposition}
\crefname{proposition}{Proposition}{Propositions}

\theoremstyle{defstyle}
\newtheorem{assumption}{Assumption}
\theoremstyle{remark}
\newtheorem{remark}{Remark}

\newcommand{\abs}[1]{\left|#1 \right|}
\newcommand{\brk}[1]{\left[#1 \right]}
\newcommand{\cB}[0]{\mathcal{B}}
\newcommand{\ceil}[1]{\left\lceil #1 \right\rceil}
\newcommand{\cF}[0]{\mathcal{F}}
\newcommand{\cR}[0]{\mathcal{R}}
\newcommand{\crl}[1]{\left\{#1 \right\}}
\newcommand{\cX}[0]{\mathcal{X}}
\newcommand{\cY}[0]{\mathcal{Y}}
\DeclareMathOperator*{\EE}{\mathbb{E}}
\newcommand{\fhat}[0]{\widehat{f}}
\newcommand{\NN}[0]{\mathbb{N}}
\DeclareMathOperator*{\PP}{\mathbb{P}}
\newcommand{\RR}[0]{\mathbb{R}}
\newcommand{\prn}[1]{\left(#1 \right)}
\newcommand{\Rhat}[0]{\widehat{\cR}}
\newcommand{\rS}[0]{\mathbf{S}}
\newcommand{\yhat}[0]{\widehat{y}}
\newcommand{\thetahat}[0]{\widehat{\theta}}

\makeatletter
\newcommand{\alphacmd@factory}[1]{}
\newcounter{alphacmdcounter}
\newcommand{\GenerateAlphabetCmds}[2]{%
    \renewcommand{\alphacmd@factory}[1]{%
        \expandafter\providecommand\csname #1##1\endcsname{{#2{##1}}}%
    }
    \setcounter{alphacmdcounter}{0}
    \loop
        \stepcounter{alphacmdcounter}
        \edef\alphacmd@ID{\@Alph\c@alphacmdcounter}
        \expandafter\alphacmd@factory\alphacmd@ID
    \ifnum\thealphacmdcounter<26
    \repeat
}
\newcommand{\GenerateAlphabetCmdsLower}[2]{%
    \renewcommand{\alphacmd@factory}[1]{%
        \expandafter\providecommand\csname #1##1\endcsname{{#2{##1}}}%
    }
    \setcounter{alphacmdcounter}{0}
    \loop
        \stepcounter{alphacmdcounter}
        \edef\alphacmd@ID{\@alph\c@alphacmdcounter}
        \expandafter\alphacmd@factory\alphacmd@ID
    \ifnum\thealphacmdcounter<26
    \repeat
}
\makeatother

\GenerateAlphabetCmds{c}{\mathcal}
\GenerateAlphabetCmdsLower{c}{\mathcal}

\GenerateAlphabetCmds{r}{\mathbf}
\GenerateAlphabetCmdsLower{r}{\mathbf}

\makeatletter
\newcommand\blfootnote[1]{%
  \begingroup
  \def\@thefnmark{}
  \renewcommand\@makefnmark{}
  \@footnotetext{#1}%
  \endgroup
}
\makeatother

\usepackage{geometry}
\title{Aggregation with Exponential Weights is Optimal in Expectation }

\author[1]{Mikael M{\o}ller H{\o}gsgaard}
\author[1]{Patrick Rebeschini}
\author[2]{Tobias Wegel}
\affil[1]{Department of Statistics, University of Oxford}
\affil[2]{Department of Computer Science, ETH Zurich}
\date{}

\begin{document}
\maketitle
\blfootnote{Authors are listed alphabetically.}

\begin{abstract}
    \noindent The aggregation with exponential weights (AEW) estimator is not fully understood in the basic setting of model selection aggregation with squared loss.
    In particular, whether it is minimax-rate optimal in expectation for large enough fixed temperatures and under random design has been an open problem since its introduction, which was explicitly posed by \citet{lecue2013optimality}. In this paper, we settle this problem by showing that \emph{without} requiring a Bernstein-type assumption, the AEW indeed achieves the excess risk $T \log (M) / (n+1)$ in expectation, whenever the temperature $T$ satisfies $(L^2/T)\exp(B/T)\leq \mu /2$. Here, the number of dictionary elements is $M$, the estimator has observed $n$ i.i.d.\ samples from any distribution, and the loss is assumed to be bounded by $B$, $L$-Lipschitz continuous and $\mu$-strongly convex. For squared loss, we show that $T\geq 4 b^2$ suffices when the predictions and labels are $[0,b]$-valued. Because AEW is known to be suboptimal in expectation for temperatures below some constant, this shows that AEW has a sharp phase transition when the temperature is large enough but constant, as conjectured by Lecu\'{e} and Mendelson. 
\end{abstract}

\section{Introduction and Main Results}

Let $(\cX,\Sigma)$ be some abstract measurable space, $\cY\subset \RR$ convex, and $P$ be an arbitrary distribution on $(\cX \times \cY, \Sigma \otimes \cB(\cY))$. Denote $\rS=\prn{X_i,Y_i}_{i=1}^n$ with $n\in\NN$ an i.i.d.\ sample from $P$ and let $\cF=\crl{f_1,\ldots,f_M}$ with $M\in\NN$ be an arbitrary but fixed finite dictionary of measurable functions $f_k:\cX\to\cY$. Given a loss function $\ell:\cY^2\to[0,\infty)$ where $\ell(\yhat,y)$ measures the loss of predicting $\yhat$ when the label is $y$, the \emph{model selection aggregation} problem is to learn a function $f:\cX\to \cY$ using the sample $\rS$ that achieves small risk
\begin{equation*}
    \cR_P(f):=\EE_{(X,Y)\sim P} \brk{\ell(f(X),Y)}
\end{equation*}
when compared to the best function in $\cF$, that is, to achieve small excess risk $ \cR_P(f)-\min_{1\leq k \leq M} \cR_P(f_k)$. 
In this paper, we consider loss functions that satisfy the following assumption:
\begin{assumption}
\label{asm:loss}
    Let $B,L,\mu>0$. The loss function $ \ell:\cY^2\to[0,\infty) $, satisfies that:
    \begin{enumerate}
        \item $\ell$ is bounded by $B$: for all $\yhat,y\in\cY$, $\ell(\yhat,y)\leq B$.
        \item $ \ell(\cdot,y) $ is $L$-Lipschitz continuous for any $ y\in \cY$: for all $y, \yhat, \yhat' \in \cY$, $\abs{\ell(\yhat,y)-\ell(\yhat',y)} \leq L \abs{\yhat-\yhat'}$.
        \item $ \ell(\cdot,y) $ is $\mu$-strongly convex for any $ y\in \cY $: for all $y\in \cY$, the map $\yhat\mapsto \ell(\yhat,y)-\frac{\mu}{2}\yhat^2$ is convex on $\cY$.
    \end{enumerate}
\end{assumption}
An important instance of such a loss function is the squared loss $\ell(\yhat,y)=(\yhat-y)^2$ on the label space $\cY=[0,b]$, in which case \cref{asm:loss} holds with $B=b^2$, $L=2b$ and $\mu=2$.

Model selection aggregation is well-studied \citep{audibert2007progressive,audibert2009fast,lecue2009aggregation,lecue2014optimal}, and the minimax rate in expectation for squared loss was proven by \citet{tsybakov2003optimal} to be, up to constant factors independent of $n$ and $M$,
\begin{equation*}
    \inf_{\fhat} \sup_{P, \abs{\cF}=M} \brk{\EE_{\rS\sim P^n}\brk{\cR_P(\fhat_{\rS})}-\min_{1\leq k \leq M} \cR_P(f_k)} \asymp \min\crl{1,\frac{\log(M)}{n}}.
\end{equation*}
We refer to \citet{mourtada2023local} for a recent overview of the literature on model selection aggregation.

Perhaps one of the most well-known estimators that one can apply to this setting is the PAC-Bayesian \emph{aggregation with exponential weights} (AEW) algorithm, which has its origins in online learning \citep{vovk1990aggregating,littlestone1994weighted,hoeven2018many} and has been studied extensively over the years \citep{yang2000mixing,catoni2004statistical,leung2006information,dalalyan2008aggregation,juditsky2008learning,rigollet2012sparse,alquier2021user,mourtada2023local}. Denoting the empirical risk with respect to the sample $\rS$ as $\Rhat_\rS(f) = \frac{1}{n}\sum_{i=1}^n \ell(f(X_i),Y_i)$, the AEW is defined as
\begin{equation}
\label{eqn:EW-definition}
    \fhat_T = \sum_{k=1}^M \thetahat_k f_k \qquad \text{where} \qquad \thetahat_k = \frac{\exp\big(-\tfrac{n}{T}\Rhat_\rS(f_k)\big)}{\sum_{j=1}^M\exp\big(-\tfrac{n}{T}\Rhat_\rS(f_j)\big)} \qquad \text{for all }k\in \crl{1,\ldots,M}.
\end{equation}
The AEW assigns a weight to each dictionary element which scales exponentially in the negative of the empirical risk $\Rhat_\rS$ with respect to the sample $\rS$. Here $T>0$ is a hyperparameter chosen by the estimator, and is called the \emph{temperature} of the exponential weights. The name ``temperature'' has its origin in the fact that the exponential weights can be viewed as a Gibbs posterior distribution over the dictionary, inspired by thermodynamics \citep{catoni2004statistical}. The larger the temperature, the more uniform the weights $\thetahat_k$ are. In particular, for small temperatures, the AEW estimator behaves similarly to empirical risk minimization on the dictionary, whereas larger temperatures induce a \emph{hedging} effect.

The known bounds for AEW are summarized by \citet{lecue2013optimality}; we restate them here for simplicity for squared loss on $[0,1]$. In short, it is known \citep{catoni2004statistical,audibert2007progressive,mourtada2023local} that if the AEW estimator is averaged into the \emph{progressive mixture rule} $\smash{\fhat_{\operatorname{pm}} = \frac{1}{n+1}\sum_{i=0}^n \fhat_{T}^{(i)}}$, where each $\fhat_{T}^{(i)}$ is an AEW estimator on the first $i$ samples and with large enough constant temperature $T$ (e.g., $T=8$), then the minimax rate is achieved, that is, $\EE[\cR_P(\fhat_{\operatorname{pm}})]\leq \min_{1\leq j\leq M} \cR_P(f_j) + 8\log(M)/(n+1)$.
Moreover, it is known that in \emph{fixed design} and for certain assumptions on the noise, an analogous guarantee can be obtained \emph{in-sample} by AEW with large enough fixed temperature \citep{dalalyan2008aggregation,dai2012deviation}. Finally, for random design, it has been shown that AEW can achieve the guarantee $\EE[\cR_P(\fhat_T)]\leq \min_{1\leq j\leq M} \cR_P(f_j) + C\log(M)/(n+1)$ under a \emph{Bernstein condition} (that is, assuming a favorable position of the dictionary relative to the distribution), where now $C$ depends on that assumption \citep{catoni2007pac,lecue2013optimality,alquier2021user}. The proof techniques appearing in those results are substantially different from those used in this work.

For random design and without assuming any Bernstein-type condition, only the following is known for the AEW estimator itself, as summarized in \citet{lecue2013optimality} (again for squared loss for simplicity).
For low temperatures $T\leq c_1$, where $c_1$ is a small enough constant independent of $M,n$, the AEW estimator is suboptimal both in expectation and in probability.
And for any temperature $T \leq c_2 \sqrt{n}/\log(n)$ (including moderately large temperatures), the AEW estimator is suboptimal on an event of constant probability.
Intuitively, the negative result for low temperatures is in line with the AEW mimicking empirical risk minimization, as the latter is also known to be suboptimal \citep{juditsky2008learning}.
However, to the best of our knowledge, the optimality of exponential weights in expectation and under random design when $T\geq c_3$ for some constant $c_3$ has remained an open problem, posed explicitly by \citet{lecue2013optimality}, who also conjectured the answer to be positive:\footnote{A previous version of the paper \citet{lecue2013optimality} states the open question of optimality as ``Question 1.2'' and explicitly conjectures the phase transition proved in the present work, as well as providing some more commentary. This previous version can be found at \href{https://maths-people.anu.edu.au/~mendelso/papers/LM6-07-07-10.pdf}{https://maths-people.anu.edu.au/\%7Emendelso/papers/LM6-07-07-10.pdf}.}

\begin{quote}
\textbf{Open Question:} \emph{In random design model selection aggregation with squared loss, is there a universal constant $c_3>0$ such that for $T\geq c_3$, the AEW estimator with temperature $T$ achieves the optimal rate of aggregation $\log(M)/n$ in expectation, uniformly over dictionaries $\cF$ of size $M$ and arbitrary distributions $P$ on $\cX\times [0,b]$, without a Bernstein condition?}
\end{quote}

In this paper, we give a positive resolution to the question, showing that for sufficiently large constant temperatures, the exponential weights estimator is optimal in expectation. We prove this under the more general \cref{asm:loss}, which contains the squared loss as a special case; however, for squared loss, we also provide a more direct proof that yields a slightly tighter bound (in terms of constant factors).

\newpage
\begin{theorem}[Minimax-rate optimality of AEW in expectation]
\label{thm:EW-optimality-expectation}
    Let the loss $\ell$ satisfy \cref{asm:loss} with parameters $B,L$ and $\mu$.
    For every $M,n\in\NN$, temperature $T\in(0,\infty)$ satisfying $(L^2/T)\exp(B/T)\leq \mu/2$, dictionary $\cF=\crl{f_1,\ldots,f_M}$ of measurable functions $f_k:\cX\to \cY$, and distribution $P$ on $\cX\times \cY$, the aggregation with exponential weights estimator \eqref{eqn:EW-definition} satisfies
    \begin{equation*}
        \EE_{\rS \sim P^n} \brk{\cR_P(\fhat_T)} \leq \min_{1\leq k \leq M} \cR_P(f_k) +\frac{T\log M}{n+1}.
    \end{equation*}
    For squared loss on $\cY=[0,b]$, the condition $(L^2/T)\exp(B/T)\leq \mu/2$ can be replaced by $T\in[4b^2,\infty)$.
\end{theorem}
\vspace{0.2cm}
\begin{remark}
    For squared loss on $[0,b]$, the parameters from \cref{asm:loss} are given by $B=b^2$, $L=2b$ and $\mu =2$, and so the first condition of \cref{thm:EW-optimality-expectation} yields $T\geq b^2/W(1/4) \approx 4.904 b^2$ where $W$ is the Lambert-W-function; our proof in the special case of squared loss hence lowers that bound on $T$ to be just $4b^2$.
    In particular, for squared loss on $[0,1]$, choosing $T=4$ yields the excess risk $4\log(M)/(n+1)$.
\end{remark}

Note that in \cref{thm:EW-optimality-expectation} no Bernstein condition is assumed. This resolves the open question from above.

\begin{figure}
    \centering
    \begin{tikzpicture}[x=1.35cm,y=1.25cm,>=stealth,font=\small]

\definecolor{subopt}{RGB}{200,60,60}
\definecolor{opt}{RGB}{60,150,80}
\definecolor{unknown}{RGB}{180,180,180}

\def\xzero{0}
\def\xcOne{3.2}
\def\xcThree{4}
\def\xgrow{7}
\def\xprob{7.5}
\def\xend{10}

\def\yexp{1.4}
\def\yprob{0}


\node[left] at (\xzero-0.2,\yexp) {in expectation};
\node[left] at (\xzero-0.2,\yprob) {in probability};

\draw[gray!50, line width=1pt] (\xzero,\yexp) -- (\xend,\yexp);

\draw[subopt, line width=5pt] (\xzero,\yexp) -- (\xcOne,\yexp);
\node[subopt, above=4pt, align=center] at ({(\xzero+\xcOne)/2},\yexp) {suboptimal\\\citep{lecue2013optimality}};
\node[below=4pt, align=center] at ({(\xzero+\xcOne)/2},\yexp)
{$T\leq c_1$};

\draw[unknown, dashed, line width=4pt] (\xcOne,\yexp) -- (\xcThree,\yexp);
\node[unknown, above=4pt] at ({(\xcOne+\xcThree)/2},\yexp) {gap};

\draw[opt, line width=5pt] (\xcThree,\yexp) -- (\xgrow,\yexp);
\node[opt, above=4pt, align=center] at ({(\xcThree+\xgrow)/2},\yexp) {optimal\\(\cref{thm:EW-optimality-expectation})};
\node[below=4pt, align=center] at ({(\xcThree+\xgrow)/2},\yexp)
{$T\geq c_3$ and constant};

\draw[subopt, line width=5pt] (\xgrow,\yexp) -- (\xend,\yexp);
\node[subopt, above=4pt, align=center] at ({(\xgrow+\xend)/2},\yexp) {suboptimal\\(\cref{prop:EW-suboptimal-growing-temperature})};
\node[below=4pt, align=center] at ({(\xgrow+\xend)/2},\yexp)
{$T\to\infty$ as $n\to\infty$};


\draw[gray!50, line width=1pt] (\xzero,\yprob) -- (\xend,\yprob);

\draw[subopt, line width=5pt] (\xzero,\yprob) -- (\xprob,\yprob);
\node[subopt, above=4pt, align=center] at ({(\xzero+\xprob)/2},\yprob) {suboptimal\\\citep{lecue2013optimality}};

\draw[subopt, line width=5pt] (\xprob,\yprob) -- (\xend,\yprob);
\node[subopt, above=4pt, align=center] at ({(\xprob+\xend)/2},\yprob) {suboptimal\\(\cref{prop:EW-suboptimal-super-high-temperatures})};

\node[below=4pt] at ({(\xzero+\xprob)/2},\yprob) {$T \le c_2\sqrt{n}/\log n$};
\node[below=4pt] at ({(\xprob+\xend)/2},\yprob) {$T \ge c_2\sqrt{n}/\log n$};

\foreach \x/\lab in {
  \xcOne/$c_1$,
  \xcThree/$c_3$,
  \xgrow/ $ $
}{
  \draw[thick] (\x,\yexp-0.2) -- (\x,\yexp+0.2);
  \node[below=2pt] at (\x,\yexp-0.2) {\lab};
}

\foreach \x/\lab in {
  \xprob/$\frac{c_2\sqrt{n}}{\log(n)}$
}{
  \draw[thick] (\x,\yprob-0.2) -- (\x,\yprob+0.2);
  \node[below=2pt] at (\x,\yprob-0.2) {\lab};
}

\end{tikzpicture}
    \vspace{-0.8cm}
    \caption{Minimax rate optimality and suboptimality of the AEW estimator for squared loss as a function of the temperature $T$, when considered uniformly over $M$, dictionaries of size $M$, and distributions.}
    \label{fig:EW-optimality-overview}
\end{figure}
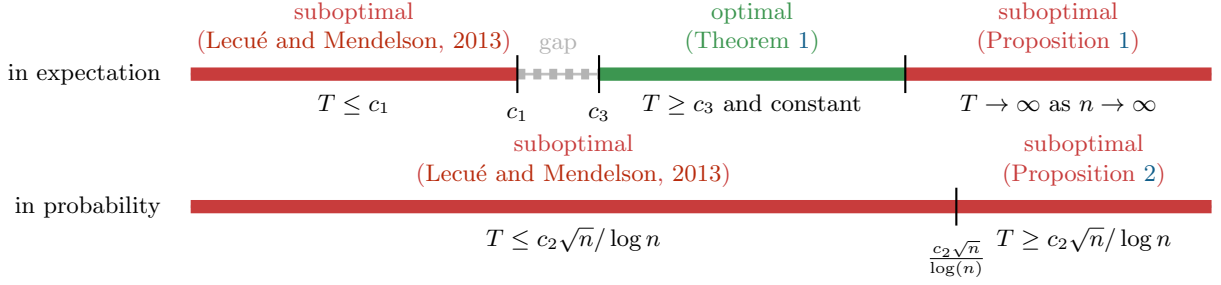
We now complement the positive result of \cref{thm:EW-optimality-expectation} by showing that if the temperature grows unboundedly with $n$ (i.e., $T\to \infty$ as $n\to \infty$) then AEW is suboptimal in expectation. See also \cref{fig:EW-optimality-overview}.
In contrast to \cref{thm:EW-optimality-expectation}, the following propositions are somewhat straightforward to prove.
\begin{proposition}
\label{prop:EW-suboptimal-growing-temperature}
    Let $\ell$ be the squared loss and $\cY=[0,1]$.
    For every $T>0$, $n\in\NN$, and $M\in\NN$ with $M\geq 2$, there exists a dictionary $\cF$ of size $M$, a distribution $P$ on $\cX\times [0,1]$, and a $\gamma_M\in[1/4,1]$ depending only on $M$, such that AEW with any temperature $T$ has excess risk lower bounded almost surely as
    \begin{equation*}
        \PP_{\rS\sim P^n}\prn{\cR_P(\fhat_{T})\geq  \min_{1\leq j\leq M} \cR_P(f_j)+ \gamma_M\min\crl{1,\frac{T\log(M-1)}{n}}}=1.
    \end{equation*}
    Moreover, it holds that $\gamma_M\to 1$ as $M\to \infty$. Now let $T_{n,M}$ be any schedule of temperatures depending on $n$ and $M$. Then, if for any sequence $M_n\geq 3$ it holds that $\log(M_n)/n\to 0$ and $T_{n,M_n} \to \infty$ as $n\to \infty$, AEW with this schedule is minimax-rate suboptimal (in expectation). Whenever $T_{n,M}$ is independent of $M$ and $T_{n,M}=T_{n}\to \infty$ as $n\to\infty$, such a sequence can be chosen so that the lower bound converges to $1$.
\end{proposition}
Along any sequence $M_n$ with $\log(M_n)/n\to 0$ and $T_{n,M_n}\to \infty$, the ratio between the lower bound in \cref{prop:EW-suboptimal-growing-temperature} and the minimax rate goes to infinity. Hence, AEW is minimax-rate suboptimal in expectation along such sequences.
\vspace{0.2cm}
\begin{remark}
    Notice that the fact that $\gamma_M\to 1$ as $M\to \infty$ implies that for large $n$ and $M$, our lower bound from \cref{prop:EW-suboptimal-growing-temperature} and the upper bound from \cref{thm:EW-optimality-expectation} witness each other's tightness not only up to constants, but in the \emph{first order}. Therefore, the universal constant factors in both bounds cannot be improved whenever $T\geq 4$. For $T\leq c_1$, the lower bound in \cref{prop:EW-suboptimal-growing-temperature} is necessarily \emph{not} tight in first order, due to the stronger lower bound by \citet{lecue2013optimality}.
\end{remark}
By a similar argument, we now complement the lower bound in probability from \citet{lecue2013optimality} for $T\leq c_2\sqrt{n}/\log(n)$ by showing that, unsurprisingly, choosing $T$ larger does not help in general. 
\begin{proposition}
\label{prop:EW-suboptimal-super-high-temperatures}
Let $\ell$ be the squared loss and $\cY=[0,1]$.
For every $n\ge 2$ and every $T\ge c_{2}\sqrt{n}/\log (n)$ (where $ c_{2} $ is the universal constant from \citet{lecue2013optimality}), there exist $M\le n^{1/c_{2}}+2$, a dictionary $\mathcal F=\{f_1,\dots,f_M\}$ with values in $[0,1]$, and a distribution $P$ on $\mathcal X\times[0,1]$ such that
\begin{equation*}
    \PP_{\rS\sim P^n}\prn{\cR_P(\fhat_T)\geq \min_{1\le j\le M}\cR_P(f_j) + \frac{1}{4\sqrt n}}=1.
\end{equation*}
\end{proposition}
Notice that since $M\leq n^{1/c_2}+2$, the dictionary size can grow at most polynomially in the sample size, and therefore the optimal rate of aggregation would be $\log(M)/n \lesssim \log(n)/n \ll 1/\sqrt{n}$.  Hence, the AEW estimator is minimax rate suboptimal for temperatures in this regime.

Together with the known results described above, this fills in the remaining gaps in the understanding of the AEW estimator for squared loss on $[0,1]$, up to the difference in constants. This is visualized in \cref{fig:EW-optimality-overview}.
Indeed, the AEW is always suboptimal in probability, and for temperatures below some constant $c_1$ or growing with $n$, it is also suboptimal in expectation. But for $T$ large enough and constant, it is optimal in expectation. This covers all values of $T$, except for the gap between these constants.

We now prove \cref{thm:EW-optimality-expectation} in \cref{proof:EW-optimality-expectation} and then \cref{prop:EW-suboptimal-growing-temperature,prop:EW-suboptimal-super-high-temperatures} in \cref{proof:EW-suboptimal-growing-temperature,proof:EW-suboptimal-super-high-temperatures}.

\section{Proof of \texorpdfstring{\cref{thm:EW-optimality-expectation}}{Theorem \ref{thm:EW-optimality-expectation}}}
\label{proof:EW-optimality-expectation}
The proof of \cref{thm:EW-optimality-expectation} reduces to a deterministic leave-one-out bound that uses ideas of \emph{algorithmic stability}. Recall that the loss $\ell$ satisfies \cref{asm:loss} with parameters $B,L,\mu$.

Let $N=n+1\geq 2$ and $\prn{x_i,y_i}_{i=1}^N\subset \cX\times \cY$ be an arbitrary fixed deterministic sample. Denote $\ell_{ij}=\ell(f_j(x_i),y_i)$, the loss of the $ i $-th observation for function $ f_j $, and $S_j = \sum_{i=1}^N \ell_{ij}$, the total loss for function $ f_j $.
We now let $p_1^{(-i)},\ldots,p_M^{(-i)}$ be the exponential weights formed on the sub-sample of size $n$ after removing the $i$th observation and with temperature $T$, that is, for each $j\in[M]$
\begin{equation*}
    p_j^{(-i)} = \frac{\exp\prn{-(S_j-\ell_{ij})/T}}{\sum_{k=1}^M \exp\prn{-(S_k-\ell_{ik})/T}} = \frac{p_j \exp\prn{\ell_{ij}/T}}{\sum_{k=1}^M p_k \exp\prn{\ell_{ik}/T}} \quad \text{where} \quad p_j = \frac{\exp\prn{-S_j/T}}{\sum_{k=1}^M \exp\prn{-S_k/T}}.
\end{equation*}
\begin{proposition}[A deterministic leave-one-out inequality]
    \label{prop:leave-one-out}
    In the setting described above, if the loss satisfies \cref{asm:loss} with parameters $B,L,\mu$, and $T$ satisfies $(L^2/T)\exp(B/T)\leq \mu/2$, it holds that
    \begin{equation}
    \label{eqn:leave-one-out}
        \frac{1}{N} \sum_{i=1}^N \ell\prn{\sum_{j=1}^M p_j^{(-i)}f_j(x_i),y_i} \leq \min_{1\leq j\leq M} \frac{1}{N}\sum_{i=1}^N\ell\prn{f_j(x_i),y_i} + \frac{T\log M}{N},
    \end{equation}
   where $ \sum_{j=1}^M p_j^{(-i)}f_j $ is the AEW estimator on the sample with the $ i $-th observation removed. For squared loss on $\cY=[0,b]$, the condition $(L^2/T)\exp(B/T)\leq \mu/2$ can be replaced by $T\in[4b^2,\infty)$.
\end{proposition}
Before proving \cref{prop:leave-one-out}, we show how it implies \cref{thm:EW-optimality-expectation}. Let $\rS'=\prn{X_i,Y_i}_{i=1}^N$ be an i.i.d.\ sample from $P$. For each $i$, we can compute the exponential weights estimator on the sample with $(X_i,Y_i)$ removed; let $\fhat_T^{(-i)}$ denote the estimator with the $i$'th observation removed and  $\thetahat^{(-i)}_j$ its weight on $f_j$. By exchangeability and independence of the data points, we  can rewrite the risk of $\fhat_T=\fhat_T^{(-N)}$ as
\begin{align*}
    \EE_{\rS\sim P^n}\brk{\cR_P(\fhat_T)}&=\EE_{\rS'\sim P^N}\brk{\ell(\fhat_T^{(-N)}(X_N),Y_N)} \\
    &= \EE_{\rS'\sim P^N}\brk{\ell(\fhat_T^{(-i)}(X_i),Y_i)} \tag{for any $i\in [N]$} \\
    &= \EE_{\rS'\sim P^N}\brk{\frac{1}{N}\sum_{i=1}^N\ell(\fhat_T^{(-i)}(X_i),Y_i)} \\
    &= \EE_{\rS'\sim P^N}\brk{\frac{1}{N}\sum_{i=1}^N\ell\prn{\sum_{j=1}^M \thetahat_j^{(-i)}f_j(X_i),Y_i}}.
\end{align*}
We can now apply \cref{prop:leave-one-out} with the sample $(X_i,Y_i)_{i=1}^N$ and $p_j^{(-i)}=\thetahat_j^{(-i)}$. By \cref{eqn:leave-one-out}, we observe that, as long as $(L^2/T)\exp(B/T)\leq \mu/2$, or $ T\geq 4b^2 $ for squared loss, the right-hand side is bounded above by
\begin{align*}
    \EE_{\rS\sim P^n}\brk{\cR_P(\fhat_T)} &\leq \EE_{\rS'\sim P^N}\brk{ \min_{1\leq j\leq M} \frac{1}{N}\sum_{i=1}^N \ell(f_j(X_i),Y_i) + \frac{T\log M}{N}} \\
    &\leq  \min_{1\leq j\leq M} \frac{1}{N}\sum_{i=1}^N \EE_{(X_i,Y_i)\sim P}\brk{\ell(f_j(X_i),Y_i)} + \frac{T\log M}{N} \\
    &= \min_{1\leq j\leq M} \cR_P(f_j) + \frac{T \log M}{n+1},
\end{align*}
which concludes the proof of \cref{thm:EW-optimality-expectation}. It remains to prove \cref{prop:leave-one-out}.

\subsection{Proof of \texorpdfstring{\cref{prop:leave-one-out}}{Proposition \ref{prop:leave-one-out}}}

Recall the notation $\ell_{ij}=\ell(f_j(x_i), y_i)$ and $S_j = \sum_{i=1}^N \ell_{ij}$.
In the first step, we apply the following lemmas. Here we make the case distinction between general losses satisfying \cref{asm:loss} and squared loss, as the proof for both is somewhat different.
They are the key to our proof and seem to be novel in the literature.

\begin{lemma}[A tilting inequality]
\label{lem:tilt-inequality-general}
    Suppose \cref{asm:loss} holds.
    Let $M\in\NN$, $(p_1,\ldots,p_M)$ be a probability distribution, and $y,\yhat_1,\ldots \yhat_M\in\cY$ be any fixed values.
    Define the \emph{tilted} probability distribution $(q_1,\ldots,q_M)$ as
    \begin{equation*}
        q_j = \frac{p_j \exp(\tfrac{1}{T} \ell(\yhat_j,y))}{\sum_{k=1}^M p_k \exp(\tfrac{1}{T} \ell(\yhat_k,y))}.
    \end{equation*}
    If $T>0$ satisfies $(L^2/T)\exp(B/T)\leq \mu/2 $, then it holds that
    \begin{equation}
    \label{eqn:tilt-inequality-general}
        \ell\prn{\sum_{j=1}^{M}q_{j}\yhat_j,y} \leq \sum_{j=1}^{M}p_{j}\ell(\yhat_j,y).
    \end{equation}
\end{lemma}
The proof of \cref{lem:tilt-inequality-general} is in \cref{proof:tilt-inequality-general}.
For squared loss we can prove the same with a slightly weaker requirement on the temperature. We split this into a separate lemma, because the proof is quite different. The proof of \cref{lem:squared-tilt-inequality} is in \cref{proof:squared-tilt-inequality}.
\begin{lemma}[A tilting inequality for squared loss]
    \label{lem:squared-tilt-inequality}
    Let $M\in\NN$, $b>0$, and $(p_1,\ldots,p_M)$ be a probability distribution, and $y,\yhat_1,\ldots,\yhat_M\in[0,b]$ be any fixed values. Define the \emph{tilted} probability distribution $(q_1,\ldots,q_M)$ as
    \begin{equation*}
        q_j = \frac{p_j \exp(\tfrac{1}{T} (\yhat_j-y)^2)}{\sum_{k=1}^M p_k \exp(\tfrac{1}{T} (\yhat_k-y)^2)}.
    \end{equation*}
    If $T\in[4b^2,\infty)$, then the following inequality is true:
    \begin{equation}
    \label{eqn:squared-tilt-inequality}
        \prn{\sum_{j=1}^M q_j \yhat_j-y}^2 \leq \sum_{j=1}^M p_j (\yhat_j-y)^2.
    \end{equation}
\end{lemma}

We use \cref{lem:tilt-inequality-general,lem:squared-tilt-inequality} with $\yhat_j=f_j(x_i),y=y_i$ by recalling that the weights of the exponential weights estimator with the $i$th datapoint removed are given by
\begin{equation*}
    p_j^{(-i)} = \frac{\exp\prn{-(S_j-\ell_{ij})/T}}{\sum_{k=1}^M \exp\prn{-(S_k-\ell_{ik})/T}} = \frac{p_j \exp\prn{ \ell_{ij}/T}}{\sum_{k=1}^M p_k \exp\prn{ \ell_{ik}/T}} = \frac{p_j \exp\prn{\ell(f_j(x_i),y_i)/T}}{\sum_{k=1}^M p_k \exp\prn{\ell(f_k(x_i),y_i)/T}},
\end{equation*}
where in turn $p_j$ are the weights of the exponential weights estimator with the full sample.
Since we assumed the necessary conditions on the temperature $T$, for a single datapoint $i$ we obtain from \cref{lem:tilt-inequality-general,lem:squared-tilt-inequality} that
\begin{align*}
    \ell\prn{\sum_{j=1}^M p_j^{(-i)}f_j(x_i),y_i} 
    \leq \sum_{j=1}^M p_j\ell(f_j(x_i),y_i) = \sum_{j=1}^M p_j\ell_{ij}.
\end{align*}
We can sum this inequality over samples $i$ and obtain that
\begin{equation*}
    \sum_{i=1}^N \ell\prn{\sum_{j=1}^M p_j^{(-i)}f_j(x_i),y_i} \leq \sum_{j=1}^M p_j \sum_{i=1}^N\ell_{ij} = \sum_{j=1}^M p_j S_j.
\end{equation*}
To conclude the proof, we can now apply the following elementary lemma to the right-hand side. It is well-known; we restate and prove (in \cref{proof:variational-exponential}) it for completeness.
\begin{lemma}[An elementary variational inequality for exponential weights]
\label{lem:variational-exponential}
    Let $M\in\NN$, let $S_1,\ldots, S_M\in \RR$ be real numbers, let $T>0$ and define $p_j = \frac{\exp\prn{-S_j/T}}{\sum_{k=1}^M \exp\prn{-S_k/T}}$ for $j=1,\ldots, M$. Then it holds that
    \begin{equation*}
        \sum_{j=1}^M p_j S_j \leq \min_{1\leq j \leq M} S_j + T\log M.
    \end{equation*}
\end{lemma}

Applying \cref{lem:variational-exponential}, to the last display and dividing by $N$ yields \cref{prop:leave-one-out}.

\subsection{Proof of \texorpdfstring{\cref{lem:tilt-inequality-general}}{Lemma \ref{lem:tilt-inequality-general}}}
\label{proof:tilt-inequality-general}

To start, we rederive the following well-known fact: due to the strong convexity from \cref{asm:loss} and convexity of $\cY$, it holds that
\begin{align}
    \ell\prn{\sum_{j=1}^M p_j \yhat_j,y} &= \brk{\ell\prn{\sum_{j=1}^M p_j \yhat_j,y}
    -\frac{\mu}{2}\prn{\sum_{j=1}^M p_j \yhat_j}^2}
    +\frac{\mu}{2}\prn{\sum_{j=1}^M p_j \yhat_j}^2 \nonumber \\
    &\leq
    \sum_{j=1}^M p_j
    \brk{\ell(\yhat_j,y)-\frac{\mu}{2}\yhat_j^2}
    +\frac{\mu}{2}\prn{\sum_{j=1}^M p_j \yhat_j}^2  \tag{Jensen's inequality and strong convexity}\\
    &=
    \sum_{j=1}^M p_j \ell(\yhat_j,y)
    -\frac{\mu}{2}
    \brk{
    \sum_{j=1}^M p_j \yhat_j^2
    -
    \prn{\sum_{j=1}^M p_j \yhat_j}^2
    } \nonumber \\
    &= \sum_{j=1}^M p_j \ell(\yhat_j,y)-\frac{\mu}{2} \sum_{j=1}^M p_j \prn{\yhat_j-\sum_{k=1}^M p_k\yhat_k}^2. \label{eqn:Jensen-gap-strong-convexity}
\end{align}
Here the last step uses the fact that for any real-valued random variable $X$, it holds $\EE[(X-\EE [X])^2]=\EE [X^2]-(\EE [X])^2$.
See also \citet[Proposition 2]{lecue2014optimal} for an analogous statement.

We can then apply Lipschitz continuity and strong convexity (through \eqref{eqn:Jensen-gap-strong-convexity}) from \cref{asm:loss} as follows:
\begin{align}
 \ell\prn{\sum_{j=1}^{M}q_{j}\yhat_j,y}
 &=\ell\prn{\sum_{j=1}^{M}q_{j}\yhat_j,y}
 -\ell\prn{\sum_{j=1}^{M}p_{j}\yhat_j,y}
 +\ell\prn{\sum_{j=1}^{M}p_{j}\yhat_j,y}
 \nonumber \\
&\leq L\abs{\sum_{j=1}^{M}q_{j}\yhat_j-\sum_{j=1}^{M}p_{j}\yhat_j}
    +\sum_{j=1}^{M}p_{j}\ell(\yhat_j,y)-\frac{\mu}{2}\sum_{j=1}^{M}p_{j}\prn{\yhat_j-\sum_{k=1}^{M}p_{k}\yhat_k}^{2}
    \tag{By $L$-Lipschitz continuity and \eqref{eqn:Jensen-gap-strong-convexity}}
\nonumber \\
&= L\abs{\sum_{j=1}^{M}q_{j}\left(\yhat_j-\sum_{k=1}^{M}p_{k}\yhat_k\right)}
    +\sum_{j=1}^{M}p_{j}\ell(\yhat_j,y)-\frac{\mu}{2}\sum_{j=1}^{M}p_{j}\prn{\yhat_j-\sum_{k=1}^{M}p_{k}\yhat_k}^{2}. \label{eqn:loss-decomposition}
\end{align}
Let now $ Z= \sum_{j=1}^{M} p_{j}\exp{(\frac{1}{T}\ell(\yhat_j,y)) } $ be the normalization constant of the distribution $(q_1,\ldots,q_M)$. To bound the first term in \eqref{eqn:loss-decomposition}, we can use
\begin{align}\label{eqn:meanvaluetheorystep}
 \sum_{j=1}^{M}q_{j}\left(\yhat_j-\sum_{k=1}^{M}p_{k}\yhat_k\right)
 &=\sum_{j=1}^{M}\frac{ \exp{(\frac{1}{T}\ell(\yhat_j,y)) } }{Z}p_{j}\left(\yhat_j-\sum_{k=1}^{M}p_{k}\yhat_k\right)\nonumber
 \\
 &=\frac{1}{Z}\sum_{j=1}^{M}p_{j}\brk{\exp\prn{\frac{1}{T}\ell(\yhat_j,y)}-\exp\prn{\frac{1}{T}\ell\prn{\sum_{k=1}^{M}p_{k}\yhat_k,y} }}\left(\yhat_j-\sum_{k=1}^{M}p_{k}\yhat_k\right)
\end{align}
where the last equality follows by $ \sum_{j=1}^{M} c p_{j}(\yhat_j-\sum_{k=1}^{M}p_{k}\yhat_k)=0 $ for any $ c\in \mathbb{R} $, especially for $ c=\exp{(\frac{1}{T}\ell(\sum_{k=1}^{M}p_{k}\yhat_k,y)) } $. Taking the absolute value on both sides and using the triangle inequality, we obtain that the first term of \eqref{eqn:loss-decomposition} is bounded by
\begin{align*}
    \abs{\sum_{j=1}^{M}q_{j}\left(\yhat_j-\sum_{k=1}^{M}p_{k}\yhat_k\right)}
    \leq\sum_{j=1}^{M}\frac{p_{j}}{Z}\abs{\exp\prn{\frac{1}{T}\ell(\yhat_j,y)}-\exp\prn{\frac{1}{T}\ell\prn{\sum_{k=1}^{M}p_{k}\yhat_k,y} }}\abs{\yhat_j-\sum_{k=1}^{M}p_{k}\yhat_k}.
\end{align*}
One can verify that the Mean Value Theorem implies the inequality $$ \forall a,b\in \RR, \forall T>0:\qquad \abs{\exp\prn{\frac{a}{T}}-\exp\prn{\frac{b}{T}}}=\abs{\int_b^a\frac{1}{T}\exp\prn{\frac{x}{T}}dx}\leq \frac{1}{T} \exp\prn{\frac{\max\crl{a,b}}{T} } |a-b|. $$ Applying this to \cref{eqn:meanvaluetheorystep}, we further obtain using $B$-boundedness of the loss that
\begin{align*}
  &\abs{\sum_{j=1}^{M}q_{j}\left(\yhat_j-\sum_{k=1}^{M}p_{k}\yhat_k\right)}\\
 &\leq \sum_{j=1}^{M}\frac{p_{j}}{TZ} \exp\prn{\frac{1}{T}\max\crl{\ell(\yhat_j,y),\ell\prn{\sum_{k=1}^{M}p_{k}\yhat_k,y} }} \cdot\abs{\ell(\yhat_j,y)-\ell\prn{\sum_{k=1}^{M}p_{k}\yhat_k,y} }\cdot\abs{\yhat_j-\sum_{k=1}^{M}p_{k}\yhat_k}
 \\
&\leq \frac{L}{T} \exp\prn{\frac{B}{T}} \sum_{j=1}^{M}\frac{p_{j}}{Z}   \prn{\yhat_j-\sum_{k=1}^{M}p_{k}\yhat_k}^{2} \tag{By $L$-Lipschitz continuity of $ \ell $, and $ \ell\leq B $}
 \\
 &\leq \frac{L}{T} \exp\prn{\frac{B}{T}} \sum_{j=1}^{M} p_{j}   \prn{\yhat_j-\sum_{k=1}^{M}p_{k}\yhat_k}^{2}
\end{align*}
where the last inequality holds because $Z=\sum_{j=1}^{M} \exp{(\frac{1}{T}\ell(\yhat_j,y)) } p_{j} \geq \sum_{j=1}^M p_j= 1$.
Thus, plugging this bound into the display \eqref{eqn:loss-decomposition}, we obtain that
\begin{align*}
     \ell\prn{\sum_{j=1}^{M}q_{j}\yhat_j,y}
     &\leq L\abs{\sum_{j=1}^{M}q_{j}\left(\yhat_j-\sum_{k=1}^{M}p_{k}\yhat_k\right)}
    +\sum_{j=1}^{M}p_{j}\ell(\yhat_j,y)-\frac{\mu}{2}\sum_{j=1}^{M}p_{j}\prn{\yhat_j-\sum_{k=1}^{M}p_{k}\yhat_k}^{2}
    \\
    &\leq \prn{\frac{L^{2}}{T}\exp\prn{\frac{B}{T}} -\frac{\mu}{2}} \sum_{j=1}^{M} p_{j}   \prn{\yhat_j-\sum_{k=1}^{M}p_{k}\yhat_k}^{2}
    +\sum_{j=1}^{M}p_{j}\ell(\yhat_j,y)
\end{align*}
where the condition on $ T $ implies $ \frac{L^{2}}{T}\exp{(\frac{B}{T})} -\frac{\mu}{2} \leq 0 $, so the last display in the above is bounded by the last term, which concludes the proof of \cref{lem:tilt-inequality-general}.

\subsection{Proof of \texorpdfstring{\cref{lem:squared-tilt-inequality}}{Lemma \ref{lem:squared-tilt-inequality}}}
\label{proof:squared-tilt-inequality}
Let $U$ be a random variable that takes the value $u_j:= \yhat_j-y \in[-b,b]$ with probability $p_j$ and denote $V=\abs{U}$ and $r=(\EE V^2)^{1/2}$, both taking values in $[0,b]$.
If $r=0$, then since $V^2=U^2\geq 0$ is a non-negative random variable we know that $U=0$ $p$-almost surely, that is, $ p_j u_j = 0$ for all $j\in [M]$. Therefore, we have that
\begin{equation*}
    \sum_{j=1}^M q_j u_j = \sum_{j=1}^M \frac{u_j p_j \exp\prn{u_j^2/T}}{\sum_{k=1}^M p_k \exp\prn{u_k^2/T}}  = 0,
\end{equation*}
and so the left hand side of \eqref{eqn:squared-tilt-inequality} vanishes. As the right hand side is always non-negative, \eqref{eqn:squared-tilt-inequality} follows.

Assume now that $r>0$. Define the function $\phi:[0,b]\to \RR$ as $\phi(v) = \exp\prn{v^2/T}/(v+r)$. A calculation using the quotient rule yields that its derivative is given by
\begin{equation*}
    \phi'(v) = \frac{(2v \exp\prn{v^2/T}(v+r)/T-\exp\prn{v^2/T})}{(v+r)^2} = \frac{\exp\prn{v^2/T}}{(v+r)^2}\prn{2v(v+r)/T-1}.
\end{equation*}
Since the first factor of the latter display is always positive, the sign of $\phi'$ is determined by the sign of $2v(v+r)/T-1$, and since $v,r\in [0,b]$, we know that $2v(v+r)/T\leq 4b^2/T \leq 1$ where we used the assumption that $T \geq 4b^2$. Therefore, on $[0,b]$ we have $\phi'\leq 0$ and $\phi$ is non-increasing.

By making a case distinction between $v\geq r$ and $v<r$, we can show that this monotonicity implies
\begin{equation}
\label{eqn:tilted-intermediate-1}
    (v-r)\exp\prn{ v^2/T} \leq \phi(r)(v^2-r^2).
\end{equation}
\emph{Case $v\ge r$:} In this case, because $\phi(v)\leq \phi(r)$ and $v^2-r^2\ge 0$, we know that $(v^2-r^2) \phi(v)\leq (v^2-r^2) \phi(r)$, and so
\begin{equation*}
    (v-r)\exp\prn{ v^2/T} = (v-r)(v+r)\phi(v) = (v^2-r^2) \phi(v)\leq (v^2-r^2) \phi(r).
\end{equation*}
\emph{Case $v<r$:} In this case, because $\phi(v)\geq \phi(r)$ and $v^2-r^2<0$, we know that $(v^2-r^2) \phi(v)\leq (v^2-r^2) \phi(r)$, and so by the same calculation, we have that $(v-r)\exp\prn{ v^2/T}\leq (v^2-r^2) \phi(r)$.

By taking expectation over $V$ in \eqref{eqn:tilted-intermediate-1} we obtain that
\begin{equation*}
    \EE\brk{(V-r)\exp\prn{ V^2/T}} \leq \phi(r) \prn{\EE V^2-r^2} = 0
\end{equation*}
where the equality follows by definition of $r=(\EE V^2)^{1/2}$. This implies $ \EE\brk{V \exp\prn{ V^2/T}}\leq r \EE\brk{\exp\prn{ V^2/T}} $, whereby we get
\begin{align*}
    \abs{\sum_{j=1}^M q_j u_j} = \abs{\sum_{j=1}^M \frac{u_j p_j \exp\prn{ u_j^2/T}}{\sum_{k=1}^M p_k \exp\prn{ u_k^2/T}} }
    &= \abs{\frac{\EE\brk{U \exp\prn{ U^2/T}}}{\EE\brk{\exp\prn{ U^2/T}}}} \\
    &\leq  \frac{\EE\brk{\abs{U} \exp\prn{ U^2/T}}}{\EE\brk{\exp\prn{ U^2/T}}}
    = \frac{\EE\brk{V \exp\prn{ V^2/T}}}{\EE\brk{\exp\prn{ V^2/T}}}\leq r.
\end{align*}
Squaring both sides and expanding the definition of $r^2$ we get that
\begin{equation*}
     \prn{\sum_{j=1}^M q_j u_j}^2 \leq r^2 = \EE V^2 = \sum_{j=1}^M p_j u_j^2.
\end{equation*}
This concludes the proof of \cref{lem:squared-tilt-inequality} by plugging back in $u_j=\yhat_j-y$.

\subsection{Proof of \texorpdfstring{\cref{lem:variational-exponential}}{Lemma \ref{lem:variational-exponential}}}
\label{proof:variational-exponential}
Recall that for all $j\in[M]$, we defined $p_j = \frac{\exp\prn{-S_j/T}}{\sum_{k=1}^M \exp\prn{-S_k/T}}$. Since $p_j>0$, we can take the logarithm and obtain
\begin{equation*}
    \log p_j = -\frac{S_j}{T} -\log \prn{\sum_{k=1}^M \exp\prn{-S_k/T}}.
\end{equation*}
We can multiply each side by $p_j$ and sum over $j\in[M]$ to obtain that
\begin{equation*}
    \sum_{j=1}^M p_j \log p_j = -\frac{1}{T}\sum_{j=1}^M p_j S_j - \log \prn{\sum_{k=1}^M \exp\prn{-S_k/T}},
\end{equation*}
where the last term remains unchanged because it is independent of $j$ and $\sum_{j=1}^M p_j =1$. Multiplying by $T$ and rearranging this, we obtain that
\begin{align*}
    \sum_{j=1}^M p_j S_j &= -T \sum_{j=1}^M p_j \log p_j- T \log \prn{\sum_{k=1}^M \exp\prn{-S_k/T}} \\
    &\leq T \log M -T \log \prn{\exp\prn{-\min_{1\leq k \leq M} S_k/T}}\\
    &= T \log M + \min_{1\leq j \leq M } S_j,
\end{align*}
where the first inequality follows because $0\leq -\sum_{j=1}^M p_j\log p_j\leq \log M$ (as it is the entropy of a probability distribution on $M$ points), and the second holds because $\sum_{k=1}^M \exp\prn{-S_k/T} \geq \exp\prn{-\min_{1\leq k \leq M} S_k/T}$. That concludes the proof of \cref{lem:variational-exponential}.

\section{Proofs of the Lower Bounds}

\subsection{Proof of \texorpdfstring{\cref{prop:EW-suboptimal-growing-temperature}}{Proposition \ref{prop:EW-suboptimal-growing-temperature}}}
\label{proof:EW-suboptimal-growing-temperature}

To prove \cref{prop:EW-suboptimal-growing-temperature}, we construct a dictionary and a distribution.
Fix $n\in\NN$ and $M\in\NN$ with $M\geq 2$. If $M=2$, then the lower bound is trivially true as $\log(M-1)=0$, so assume without loss of generality that $M\geq 3$, so that $\log(M-1)>0$. Define for some $\alpha\in(0,1]$ to be chosen later
\begin{equation*}
    a^2 = \alpha r \qquad \text{where} \qquad r= \min\crl{1,\frac{T\log(M-1)}{n}}.
\end{equation*}
Let the distribution $P$ of $(X,Y)$ be such that $Y= 0$ almost surely, and choose the dictionary $f_1\equiv 0$, as well as $f_2=\cdots=f_M\equiv a$. Thus $f_1$ is the optimal dictionary element, whereas all others are suboptimal. 
For every sample $\rS$, almost surely, we have that
\[
\Rhat_\rS(f_1)=0,
\qquad
\Rhat_\rS(f_j)=a^2 \quad \text{for } j\ge2.
\]
Hence the total mass put on the suboptimal functions with index $j\geq 2$ is
\[
1-\thetahat_1 = \sum_{j=2}^M \frac{\exp\prn{-na^2/T}}{1+\sum_{k=2}^M \exp\prn{-na^2/T}} = \frac{(M-1)\exp\prn{-na^2/T}}
{1+(M-1)\exp\prn{-na^2/T}}.
\]
Because $r \leq T\log (M-1)/n$ by definition of $r$, we get that $na^2/T = \alpha nr/T\leq \alpha \log(M-1)$. Therefore, we obtain that $(M-1)\exp(-na^2/T) \geq (M-1)^{1-\alpha}$ and so
\begin{equation*}
    1-\thetahat_1 \geq \frac{(M-1)^{1-\alpha}}{1+(M-1)^{1-\alpha}}.
\end{equation*}
Notice that because the AEW estimator on this dictionary is $\fhat_{T} \equiv a(1-\thetahat_1)$, and because $f_1$ has risk zero, this means it has $P^n$-almost surely an excess risk of at least 
\begin{equation*}
    \cR_P(\fhat_T)- \min_{1\leq j \leq M} \cR_P(f_j) \geq a^2(1-\thetahat_1)^2 \geq \alpha \prn{\frac{(M-1)^{1-\alpha}}{1+(M-1)^{1-\alpha}}}^2 r.
\end{equation*}
Since this lower bound is true for every $\alpha\in(0,1]$, we can define 
\begin{equation*}
    \gamma_M =\sup_{\alpha\in (0,1]} \alpha \prn{\frac{(M-1)^{1-\alpha}}{1+(M-1)^{1-\alpha}}}^2
\end{equation*}
and the excess risk is lower bounded by $\gamma_M r$ for the $\alpha\in(0,1]$ attaining the maximum (which exists). Notice that $\gamma_M\in [1/4,1]$ since for $\alpha=1$ the term is $1/4$, so $\gamma_M$ must be larger, and $\gamma_M\leq 1$, follows from $\alpha\leq 1$ and $((M-1)^{1-\alpha}/(1+(M-1)^{1-\alpha})\leq 1$.
We now show that $\gamma_M\to 1$ as $M\to \infty$. To that end, denote the function to be optimized as $h_M(\alpha) = \alpha ((M-1)^{1-\alpha}/(1+(M-1)^{1-\alpha}))^2$. To prove the limit, consider the sequence $\alpha_M$ defined as $\alpha_M = 1-1/\sqrt{\log(M-1)}$. Then we get that $(M-1)^{1-\alpha_M} =\exp\prn{\sqrt{\log(M-1)}}$ and so
\begin{equation*}
    h_{M}(\alpha_M)= \prn{1-\frac{1}{\sqrt{\log(M-1)}}}\prn{\frac{\exp\prn{\sqrt{\log(M-1)}}}{1+\exp\prn{\sqrt{\log(M-1)}}}}^2 \to 1 \quad \text{as } M\to\infty.
\end{equation*}
By definition, we then have that $\gamma_M \geq h_M(\alpha_M) \to 1$ and $\gamma_M\leq 1$, implying that $\gamma_M\to 1$ as $M\to \infty$.
Plugging in the definition of $r$ yields the proof of the first two claims in \cref{prop:EW-suboptimal-growing-temperature}.

Now, let $M_n\ge3$ be such that $\log(M_n)/n\to0$ and
$T_{n,M_n}\to\infty$. By the first part, for each $n$ there is a
distribution and dictionary for which
\[
    \cR_P(\widehat f_{T_{n,M_n}})-\min_{1\leq j\leq M} \cR_P(f_j)
    \ge \gamma_{M_n}\min\crl{1,\frac{T_{n,M_n}\log(M_n-1)}{n}}.
\]
Since $M_n\ge3$, there exists a universal constant $c>0$ such that
$\log(M_n-1)\ge c\log M_n$. Therefore
\[
    \frac{\mathcal R_P(\widehat f_{T_{n,M_n}})-\min_{1\leq j \leq M} \cR_P(f_j)}{\log(M_n)/n}
    \geq \frac{1}{4} \min\crl{\frac{n}{\log M_n}, c T_{n,M_n}}\to\infty \quad \text{as } n\to \infty.
\]
Thus, the excess risk is of order larger than $\log(M_n)/n$, and hence suboptimal.

It remains to show that whenever $T_{n,M}$ is independent of $M$ and $T_{n,M}=T_{n}\to \infty$ as $n\to\infty$, there exists a sequence $M_n\to \infty$ such that $\log (M_n) / n\to 0$ but $r=1$ for all $n$. By the first part, that then yields the lower bound $\gamma_{M_n}$ converging to $1$ because $M_n\to \infty$.
To that end, consider the sequence $M_n = \ceil{\exp(\max\crl{n/T_n,\sqrt{n}})}+1 \leq 3\exp(\max\crl{n/T_n,\sqrt{n}})$. Then
\begin{equation*}
    \frac{\log(M_n)}{n} \leq \max\crl{\frac{1}{T_n},\frac{1}{\sqrt{n}}} + \frac{\log 3}{n} \to 0 \quad \text{as } n\to \infty,
\end{equation*}
and at the same time, $r=1$ because the other term in the minimum of the definition of $r$ is lower bounded by
\begin{equation*}
    \frac{T_n \log(M_n-1)}{n} \geq \frac{T_n \log(\ceil{\exp(\max\crl{n/T_n,\sqrt{n}})})}{n} \geq \frac{T_n(n/T_n)}{n} =1.
\end{equation*}
That concludes the proof of \cref{prop:EW-suboptimal-growing-temperature}.

\subsection{Proof of \texorpdfstring{\cref{prop:EW-suboptimal-super-high-temperatures}}{Proposition \ref{prop:EW-suboptimal-super-high-temperatures}}}
\label{proof:EW-suboptimal-super-high-temperatures}

To prove \cref{prop:EW-suboptimal-super-high-temperatures}, we construct a dictionary and a distribution. We proceed almost identically to \cref{proof:EW-suboptimal-growing-temperature}.
Let the distribution $P$ of $(X,Y)$ be such that $Y= 0$ almost surely, and choose the dictionary $f_1\equiv 0$, as well as $f_2=\cdots=f_M\equiv a$, where we choose $a>0$ as $a^2=n^{-1/2}$.
Moreover, we choose the number of dictionary elements $M$ to satisfy
\[
M-1=\ceil{\exp\prn{\sqrt n/T}}.
\]
The assumption on $T$ that $T\geq c_2\sqrt{n}/\log(n)$ gives 
\begin{equation*}
    M=\ceil{\exp\prn{\sqrt n/T}}+1\leq \ceil{\exp\prn{\log(n)/c_2}}+1 \le n^{1/c_{2}}+2,
\end{equation*}
confirming that the dictionary is not too large.
For every sample $\rS$, we have by construction that $P$-almost surely,
\[
\Rhat_\rS(f_1)=0,
\qquad
\Rhat_\rS(f_j)=a^2 \quad \text{for } j\ge2.
\]
Hence the total mass put on the suboptimal functions with index $j\geq 2$ is
\[
1-\thetahat_1 = \sum_{j=2}^M \frac{\exp\prn{-na^2/T}}{1+\sum_{k=2}^M \exp\prn{-na^2/T}} = \frac{(M-1)\exp\prn{-na^2/T}}
{1+(M-1)\exp\prn{-na^2/T}}
=
\frac{(M-1)\exp\prn{-\sqrt n/T}}
{1+(M-1)\exp\prn{-\sqrt n/T}}
\ge \frac{1}{2}
\]
where the last inequality follows from $M-1\geq \exp\prn{\sqrt{n}/T}$ implying $(M-1)\exp\prn{-\sqrt n/T}\geq 1$. Therefore, by noticing that $\fhat_T\equiv (1-\thetahat_1)a$, and since $f_1$ has risk zero,
\[
\cR_P(\fhat_T)-\min_{1\leq j\leq M}\cR_P(f_j)=((1-\thetahat_1)a)^2
\ge \frac{1}{4} n^{-1/2}.
\]
The above holds almost surely, so the event has probability one. That concludes the proof.
\section{Discussion}

In this short paper, we prove that the aggregation with exponential weights estimator achieves the minimax optimal rate of aggregation $T\log(M)/(n+1)$ with respect to $M$ and $n$, for large enough fixed temperatures $T$, when the loss is bounded, Lipschitz continuous, and strongly convex. Importantly, this does not require a Bernstein condition and includes the squared loss as a special case.

The proof at its core uses an average leave-one-out stability argument (\cref{prop:leave-one-out}). The reduction to such a leave-one-out bound is similar in spirit to existing bounds, see for instance \citet{ForsterW02,koren2015fast}. However, the main difference and novelty of our approach is in the proof of that stability result, specifically in \cref{lem:tilt-inequality-general,lem:squared-tilt-inequality}, which explicitly makes use of the tilting of the exponential weights distribution when one sample is left out.
\cref{lem:tilt-inequality-general,lem:squared-tilt-inequality} seem to be novel, and they may also be of independent interest.

We would like to remark that the tightness of \cref{lem:tilt-inequality-general,lem:squared-tilt-inequality} is crucial. Indeed, a weaker version of \cref{lem:tilt-inequality-general,lem:squared-tilt-inequality} can be obtained via an application of Jensen's inequality. In the notation of \cref{lem:tilt-inequality-general},
\begin{equation*}
    \ell\prn{\sum_{j=1}^M q_j \yhat_j,y}\leq \sum_{j=1}^M q_j \ell(\yhat_j,y) = \sum_{j=1}^M \frac{p_j \exp\prn{ \ell(\yhat_j,y)/T}}{\sum_{k=1}^M p_k \exp\prn{ \ell(\yhat_k,y)/T}} \ell(\yhat_j,y) \leq e^{B/T} \sum_{j=1}^M p_j \ell(\yhat_j,y),
\end{equation*}
where the first inequality is Jensen, and the second inequality uses that the loss is bounded $\ell \leq B$ and $\sum_{k=1}^M p_k \exp\prn{\ell(\yhat_k,y)/T}\geq 1$. The resulting bound is worse than \cref{eqn:tilt-inequality-general,eqn:squared-tilt-inequality} by the factor $e^{B/T}$, which for any $T\in(0,\infty)$ is strictly larger than $1$. Importantly, any factor strictly larger than $1$ means that the final bound has the same factor in front of $\min_{1\leq j\leq M} \frac{1}{N}\sum_{i=1}^N (y_i-f_j(x_i))^2$ in \cref{prop:leave-one-out}, respectively $\min_{1\leq k \leq M} \cR_P(f_k)$ in \cref{thm:EW-optimality-expectation}. Therefore, this approach would not yield a positive conclusion to the open question. This highlights the importance of \cref{lem:tilt-inequality-general,lem:squared-tilt-inequality}.
We also remark that the phenomenon of achieving fast rates at the cost of a worse comparator is common in some PAC-Bayesian analyses of AEW, see for instance \citet[Example 3.1 and surrounding discussion]{alquier2021user}.

Together with the suboptimality results by \citet{lecue2013optimality}, \cref{thm:EW-optimality-expectation,prop:EW-suboptimal-growing-temperature} show that the AEW estimator undergoes a sharp phase transition when the temperature is constant, and in particular the exact \emph{value} of that constant is crucial. This further demonstrates the sensitivity of the AEW to the temperature parameter, as argued by \citet{lecue2013optimality}.

\section*{Acknowledgements}

The authors thank Tomas Va\v{s}kevi\v{c}ius for helpful input.
Tobias Wegel was supported by SNSF Grant 204439. Mikael M{\o}ller H{\o}gsgaard was supported by a Carlsberg Internationalisation Fellowship. Patrick Rebeschini was funded by UK Research and Innovation (UKRI) under the UK government’s Horizon Europe funding guarantee [grant number EP/Y028333/1].

\paragraph{LLM Usage.}
The proof idea of \cref{thm:EW-optimality-expectation}, and specifically \cref{lem:squared-tilt-inequality}, is based on interactions the authors had with ChatGPT 5.5. While the authors take full responsibility for the contents of this work and the correctness of the proof, they acknowledge the significant impact the LLM had on this work.

\newpage
\bibliographystyle{apalike} 
{
\small
\bibliography{bibliographyEW}
}

\end{document}